\documentclass[reqno]{amsart}
\usepackage{amsmath,amsthm,amsfonts,amssymb}
\usepackage{xcolor}
\usepackage{booktabs,algorithm2e}
\usepackage{hyperref}

\usepackage{tikz}
\usetikzlibrary{chains,matrix,scopes,fit,decorations,positioning,shapes,arrows,decorations.markings,arrows.meta}
\newcommand\BBR{\mathbb{R}}
\newcommand\BBP{\mathbb{P}}
\newcommand\BBE{\mathbb{E}}
\newcommand\BBM{\mathbb{M}}

\newcommand\BX{\mathbf{X}}

\newcommand\cN{\mathcal{N}}

\newcommand\cM{\mathcal{M}}

\newcommand\cS{\mathcal{S}}

\newcommand\cL{\mathcal{L}}

\newcommand\ve{\varepsilon}

\renewcommand\P{\mathbb{P}}

\newcommand\R{\mathbb{R}}

\DeclareMathOperator{\Cov}{Cov}

\DeclareMathOperator{\Pa}{Pa}
\DeclareMathOperator{\Ch}{Ch}

\DeclareMathOperator{\CIS}{CIS}
\DeclareMathOperator{\MTP}{MTP}
\DeclareMathOperator{\PA}{PA}

\newtheorem{theorem}{Theorem}
\newtheorem{proposition}[theorem]{Proposition}
\newtheorem{lemma}[theorem]{Lemma}
\newtheorem{corollary}[theorem]{Corollary}

\theoremstyle{definition}
\newtheorem{definition}[theorem]{Definition}
\newtheorem{remark}[theorem]{Remark}
\newtheorem{example}[theorem]{Example}

\numberwithin{equation}{section}
\numberwithin{theorem}{section}

\title{Positivity in Linear Gaussian Structural Equation Models}
\author{Asad Lodhia}
\author{Jan-Christian H\"{u}tter}
\author{Caroline Uhler}
\author{Piotr Zwiernik}

\thanks{AL was supported by a fellowship from the Eric and Wendy Schmidt Center at the Broad Institute of MIT and Harvard. CU was partially supported by NCCIH/NIH (1DP2AT012345), ONR (N00014-22-1-2116), NSF (DMS-1651995), the MIT-IBM Watson AI Lab, and a Simons Investigator Award. PZ was supported by the NSERC Discovery Grant RGPIN-2023-03481.}

\begin{document}
\maketitle
\begin{abstract}

We study a notion of positivity of Gaussian directed 
acyclic graphical models corresponding to a non-negativity 
constraint on the coefficients of the associated structural equation 
model. We prove that this constraint is equivalent to the distribution being conditionally 
increasing in sequence (CIS), a well-known 
subclass of positively associated random variables. These distributions require 
knowledge of a permutation, a CIS ordering, of the nodes for which the 
constraint of non-negativity holds. 
We provide an algorithm and prove in the noise-less setting that a CIS ordering can be recovered when it exists. We extend this result to the noisy setting and provide assumptions for recovering the CIS orderings. In addition, we provide a characterization of Markov equivalence for CIS DAG models. 
Further, we show that when a 
CIS ordering is known, the corresponding class of Gaussians lies in a family of 
distributions in which maximum likelihood estimation is a convex 
problem.
\end{abstract}

\section{Introduction}

Many random systems exhibit some form of positive dependence. Examples 
in statistical physics include the  ferromagnetic Ising Model 
\cite{G70} as well as general classes of lattice gas models and 
percolation models \cite{FKG71}. In fields such as finance 
\cite{agrawal2022covariance}, psychometrics \cite{ST80,TR83,NRT83} and 
biology (see \cite{Ma09} and further discussion in \cite[Section 
1.1]{lauritzen2020locally}), positive dependence naturally arises~\cite{carter2022positive}. In recent years, there has been an increased 
interest in exploiting this and related notions of positive dependence 
in statistical modelling and in machine learning. 

This research direction has been particularly fruitful in the context 
of Gaussian and related distributions. Well studied examples of 
positive dependence in Gaussian models include: positive association 
defined by nonnegativity of all correlations \cite{pitt1982positively}, 
totally positive distributions (also known as MTP$_2$ distributions) defined by nonnegativity of all partial 
correlations \cite{SH15,LUZ19}, and mixtures of these two 
scenarios 
as discussed in \cite{lauritzen2020locally}. Various methods have been developed for covariance matrix estimation in the Gaussian setting~\cite{egilmez2017graph,agrawal2022covariance,soloff2020covariance,zhou2022covariance}. 
In applications, where the assumption of positive dependence is appropriate, these methods 
perform extremely well with no need for explicit regularization \cite{agrawal2022covariance,rossell2021dependence}.

An important problem, which motivates our work, is that none of these 
notions of positive dependence are suitable in the 
context of directed acyclic graphical models, also known as Bayesian networks\footnote{For  background on 
Bayesian networks and associated graphical models, see e.g.~\cite{La96}.}. 
For a simple example, consider a Gaussian vector $X=(X_1,X_2,X_3)$ such 
that all partial correlations are non-negative. In other words, 
expressing partial correlations in terms of marginal correlation 
coefficients $\rho_{ij}={\rm corr}(X_i,X_j)$, we require that 
$$\rho_{12}-\rho_{13}\rho_{23}\ge 0,\quad  
\rho_{13}-\rho_{12}\rho_{23}\ge 0,\quad  
\rho_{23}-\rho_{12}\rho_{13}\ge 0.$$ If $X$ is Markov to the DAG $2\to 
1\leftarrow 3$, or equivalently, if $\rho_{23}=0$, then, from the last 
inequality, we necessarily have that $\rho_{12}\rho_{13}=0$. In other 
words, a Bayesian network with a v-structure cannot have all partial 
correlations non-negative.  Given that two DAGs are Markov equivalent if 
and only if they have the same skeleton and v-structures (see 
Theorem~\ref{thm:markovequiv}), adding the MTP$_2$ constraint would severely restrict the class of Bayesian networks.

In this paper, we study a natural form of directional positive 
dependence that is suitable for Gaussian models on directed acyclic 
graphs (DAGs). We introduce the model through its representation via
linear structural equations~\cite{Pe09}. If $G$ is a DAG with $m$ 
nodes representing the Gaussian vector $X=(X_1,\ldots,X_m)$, then the 
distribution of $X$ lies in the associated DAG model if it admits the 
stochastic representation
\begin{equation}\label{eq:SEM}
X_i\;=\;\sum_{j\in \Pa(i)}\Lambda_{ij}X_j+\varepsilon_i 
\quad\mbox{for all }i=1,\ldots,m,
\end{equation}
where $\Pa(i)$ denotes the set of parents of the node $i$ in the DAG 
$G$, $\Lambda_{ij}\in \R$, and $\varepsilon_i\sim N(0,\sigma_i^2)$ are 
mutually independent. In matrix form, this can be written as 
$X=\Lambda X+\varepsilon$, where $\Lambda_{ij}=0$ unless $j\to i$ in 
$G$.
\begin{remark}
\label{rem:precfromlam}
Let $D$ be a diagonal matrix representing the covariance matrix of 
$\varepsilon$. Denoting the covariance matrix of $X$ by $\Sigma$ (we assume throughout that it is full rank), then \eqref{eq:SEM} implies
\[
\Sigma \;=\; (I 
-\Lambda)^{-1} D (I-\Lambda)^{-\top},
\]
which is equivalent to the following equality for the precision matrix $K = 
\Sigma^{-1}$:
\[
K \;=\; (I-\Lambda)^\top D^{-1} (I - \Lambda).
\]
\end{remark}

The following natural notion of positivity in Gaussian DAG models is 
the central theme of our paper. 

\begin{definition}\label{def:posDAG}
A Gaussian vector $X$ is \emph{positively DAG dependent} with respect to a DAG 
$G$ if $X$ admits the stochastic representation \eqref{eq:SEM} with all 
$\Lambda_{ij}$ nonnegative. We denote the subset of Gaussian DAG models $\BBM(G)$ 
over $G$ that satisfy this property by $\BBM_+(G)$.
\end{definition}

Other approaches have been proposed to define positive dependence on 
DAGs; see, for example, 
\cite{vanderweele2010signed,wellman1990fundamental} and references 
therein. To explain the relationship between these different notions, we note that positive DAG dependence is closely related to 
the following classical notion of positive dependence. 

\begin{definition}
\label{def:CISrvar}
A random vector $X = (X_1,\ldots, X_m)$ is \emph{conditionally increasing in
sequence} (\emph{CIS}) if for every $i\in[m]$ and every fixed $x_i \in \BBR$,
it holds that
\[
\BBP\big(\{X_i \geq x_i \} \big| (X_j = x_j)_{j < i}\big)
\]
is a non-decreasing function in $(x_1,\ldots, x_{i-1})$, when equipped 
with the coordinatewise partial order.
\end{definition}

In the context of DAGs, the papers 
\cite{vanderweele2010signed,wellman1990fundamental} investigated a 
similar notion which they called a ``weak monotonic effect'' or 
``positive influence''. If a parent $k\in \Pa(i)$ of a particular 
vertex $i$ has the property that
\[
\BBP\Big( \{ X_i \geq x_i\} \Big| \bigcap_{j\in \Pa(i)} (X_j = x_j ) 
\Big)
\]
is a non-decreasing (non-increasing) function in  $x_k$ then $k$ is said to have a weak monotonic positive effect on $i$. Notably, this 
condition can be used to infer the presence/absence of certain edges in the graph. 
To provide a specific example, consider the DAG from
\cite[Example 5]{vanderweele2010signed} with variables A denoting 
air pollution levels, E denoting antihistamine treatment, D denoting 
asthma incidence, and C denoting bronchial reactivity in the following figure.
\begin{figure}[h]
\includegraphics{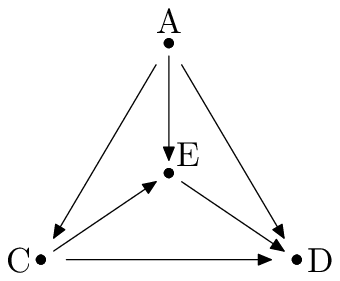}
\label{fig:weeleexamp}
\end{figure}

\noindent From context, it is reasonable to assume that the directed edges $(A,D)$, $(A,E)$ and 
$(A,C)$ are weak positive monotonic effects, and similarly the edges $(C,E)$ 
and $(C,D)$ are weak positive monotonic effects. The following argument can be used to test the causal relationship $E\rightarrow D$: 
From \cite[Theorem 4]{vanderweele2010signed} it follows that 
the covariance of $E$ and $D$ must be non-negative due to the weak 
positive monotonic effects of the other edges. Thus if the observed covariance of $E$ and $D$ was negative (which is the desired medical objective), we 
would conclude the presence of the edge $E\rightarrow D$ even without measuring the variables $A$ and $C$.

We will show that the notion of positive dependence of a DAG considered in this work
(stated in Definition~\ref{def:generalposdep}) is the same as 
assuming a weak positive monotonic effect of \emph{every} parent on its child. This example showing how positive dependence can be used to derive causal relationships motivates our study of Markov equivalence of these models. 

In this work, we link the class of CIS models 
(which a priori make no reference to any DAG structure) to positive DAG 
dependence. 
In particular, we will discuss the problem of ordering the variables in such a way that it is CIS and identifying when there exists such an ordering. The resulting DAG could be used for causal inference. 

In Theorem~\ref{th:CISinDAG} below, we show that a Gaussian vector $X$ 
is CIS if and only if it is positively DAG dependent with respect to 
the full DAG with arrows $i\to j$ for all $i < j$. It follows that the 
Gaussian CIS condition has a simple algebraic formulation. Let $K=UU^T$ 
be the Cholesky decomposition ($U$ is upper triangular with positive 
diagonal entries) and $K=\Sigma^{-1}$ is the inverse covariance matrix 
of $X$. Then our notion of positive dependence restricts the signs of 
the off-diagonal entries of $U$ to be non-positive. This constraint is 
convex, which makes computing the maximum likelihood estimator (MLE) particularly 
tractable. 

In practice, $K$ may admit such a signed Cholesky factorization only 
after permuting its rows and columns. Thus, part of the problem is to 
recover a permutation matrix $P$ that makes such a signed factorization 
possible. Maximizing the likelihood over all $m!$ permutation matrices 
is infeasible. Instead, we propose a simple algorithm for learning 
such a permutation, and we provide statistical guarantees for the proposed algorithm. 

We will often contrast the class of CIS Gaussian vectors $X$ with the 
well-known and well-studied class of \textit{multivariate totally positive distributions of 
order 2} (MTP${}_2$), which requires that its density $p$ on 
$\BBR^m$ satisfies 
\[ 
p(x) p(y) \;\leq\; p(x\vee y) p(x\wedge y) \qquad \hbox{for all $x, y 
\in \BBR^m$}, 
\] 
where $\vee$ is the componentwise maximum and $\wedge$ is the componentwise 
minimum. This inequality appeared in \cite{FKG71}, where it was shown to 
imply positive association for general distributions. In the Gaussian 
case MTP${}_2$ was shown to be equivalent to the 
precision matrix (inverse covariance matrix) being an M-matrix ~\cite{KR83}. 

\begin{definition}[M-matrix]
\label{def:m-matrix}
A positive definite $m \times m$ matrix $A= [a_{i,j}]_{1\leq i ,j \leq 
m}$ is an M-matrix if the entries satisfy $a_{i,j} \leq 0$ for all $i 
\neq j$. The space of symmetric, positive definite M-matrices of 
dimension $m \times m$ is denoted $\cM_m(\BBR)$.
\end{definition} 

\subsection{Outline} Section~\ref{sec:cisdefandprop} expounds upon the 
relationship between CIS distributions and DAG models while also 
providing motivating examples both in the Gaussian and 
non-Gaussian settings. Section~\ref{sec:examples} provides examples 
that distinguish CIS distributions from MTP${}_2$ and other positively 
associated distributions along with an illustration that CIS orderings 
may not provide sufficient information to recover the underlying Markov 
equivalence class. Section~\ref{sec:markovandcis} dives deeper into Markov equivalence for CIS 
orderings. Section~\ref{sec:choleskestimat} shifts the focus to 
parameter estimation and fitting: Cholesky factor 
models are introduced for the purpose of characterizing the MLE  of $\Lambda$ and $D$ of a CIS distributed vector assuming 
the underlying CIS ordering is known. Section~\ref{sec:cisordrecov} concerns 
recovering a CIS ordering, first in the population case and then 
proving consistency of a noisy version of our proposed algorithm under 
simple assumptions on $\Lambda$. In this section, we also prove results on what sorts of CIS orderings are possible for a 
distribution.

\subsection{Notation} For a DAG $G = (V,E)$, we denote the set of parent nodes of 
a vertex $i$ by $\Pa(i)$ and the set of children nodes 
of a vertex $i$ by $\Ch(i)$. If there are several DAGs 
over the same vertex set $V$ under consideration, we write $\Pa_G(i)$ 
and $\Ch_G(i)$ to indicate the dependence on the particular DAG $G$.
We will mostly use $V = [m] = \{1,\ldots, m\}$ or subsets of $[m]$. 

When we say a function $f: \BBR^k \to \BBR$ is increasing 
(non-decreasing) in $\BBR^k$, we mean that $f$ is increasing 
(non-decreasing) in each variable. Moreover for a subset $A \subset 
[k]$, if we write $(x_j)_{j\in A}$, or equivalently, $x_A$, we mean the 
tuple formed by taking the entries of $x$ that are indexed by $A$, 
keeping the original order.

We denote the set of $m \times m$ positive semidefinite 
matrices by $\cS_m(\BBR)$ and the subset of positive definite matrices 
by $\cS_m^+(\BBR)$. Further, $I$ always denotes the identity matrix. 
When $M$ is an $s\times t$ matrix with $A \subset [s]$ and $B\subset [t]$, 
then $M_{A,B}$ is the submatrix of size $|A| \times |B|$ with entries 
$M_{i,j}$ with $i \in A$ and $j \in B$. Following \cite[Section 
5.1.1]{La96}, if a matrix operation appears with the subset indices, 
e.g., $M^{-1}_{A,A}$ the matrix operation is performed first --- so 
$M^{-1}_{A,A}$ is the submatrix of $M^{-1}$ indexed by $A$, whereas 
$(M_{A,A})^{-1}$ is the inverse of the submatrix of $M$ indexed by $A$. 
We will use the shorthand $\backslash i$ for $[m] \backslash i$. 

When we consider collections of permutations, we use one line notation 
and use parentheses around those elements that can be ordered in any 
way, so for instance $(123)45$ is the set of permutations for which 
$\sigma(4) = 4$ and $\sigma(5) = 5$ and $1,2,3$ can be arbitrarily 
assigned to the values $\sigma(1)$, $\sigma(2)$ and $\sigma(3)$, that 
is, $(123)45=\{12345,13245,21345,23145,31245,32145\}$.

\section{Structure of Positive Dependence on a DAG}
\label{sec:cisdefandprop}

\subsection{Basic results and definitions}

We start by stating the main result of this section, which links the 
classical concept of CIS dependence and positive DAG dependence.
\begin{theorem}\label{th:CISinDAG}
A Gaussian vector $X$ is CIS if and only if it is positively DAG 
dependent with respect to the full DAG with $i\to j$ for all $i<j$.
\end{theorem}
The proof relies on a lemma that we prove first.

\begin{lemma}
\label{lem:nonnegprecis}
Let $Z =(Z_1,\ldots,Z_m) \sim \cN_m(\mu,\Sigma)$ be a Gaussian random 
vector on $\BBR^m$ with mean $\mu \in \BBR^m$ and covariance 
$\Sigma\in\cS_m^+(\BBR)$ , let $K = \Sigma^{-1}$ be the precision 
matrix. The function
\begin{equation}
\label{eq:survive}
\BBP\Big( \{Z_i \geq x_i\} \Big| \bigcap_{j \neq i} \{Z_j = x_j\} 
\Big)    
\end{equation}
is non-decreasing in $(x_j)_{j\neq i}$ if and only if  $K_{i,j} \leq 0$ 
for all $j \neq i$. Moreover, this statement is equivalent to the 
following two statements:
\begin{enumerate}
	\item [(a)] $\BBE[Z_i | Z_{\backslash i}]$ is a non-decreasing 
	function in $(Z_j)_{j\neq i}$.
	\item [(b)] $Z_i=\sum_{j\neq i} \Lambda_{ij} Z_j+\varepsilon_i$ 
	with $\Lambda_{ij}\geq 0$ and $\varepsilon_i$ Gaussian and 
	independent of $(Z_j)_{j\neq i}$. 
\end{enumerate}
\end{lemma}
\begin{proof}
 It is a classic result \cite[Theorem 1.2.11 (b)]{Mu82} that 
\[
Z_i | Z_{\backslash i} \sim \cN\Big(\mu_i + \Sigma_{i,\backslash i} 
(\Sigma_{\backslash i,\backslash i})^{-1} (Z_{\backslash i} - 
\mu_{\backslash i}), \Sigma_{i,i} - \Sigma_{i,\backslash i} 
(\Sigma_{\backslash i,\backslash i})^{-1} \Sigma_{i,\backslash i}^\top 
\Big),
\]
but note that by the Schur complement formula,
\begin{align*}
K_{i,i} &=  \Big(\Sigma_{i,i} - \Sigma_{i,\backslash i} 
\big(\Sigma_{\backslash i,\backslash i}\big)^{-1} \Sigma_{i,\backslash 
i}^\top\Big)^{-1},\\
K_{i,\backslash i} &=  - K_{i,i} \Sigma_{i,\backslash i} 
\big(\Sigma_{\backslash i,\backslash i}\big)^{-1},
\end{align*}
and $K_{i,i} > 0$ by positive definiteness.
Hence we may rewrite the mean of $Z_i | Z_{\backslash i}$ as
\[
\mu_i - \frac{K_{i,\backslash i}}{K_{i,i}} (Z_{\backslash i} - 
\mu_{\backslash{i}}).
\]
It is then clear that the function in the statement of the lemma is 
non-decreasing in $x_{\backslash i}$ only if the entries of 
$K_{i,\backslash i}$ are all non-positive. Note that this is also the 
condition on the conditional mean in (a). Equivalence with (b) follows 
from the fact that 
\[
\varepsilon_i\;:=\;Z_i\;-\;\BBE[Z_i| Z_{\backslash i}]
\]
is a mean zero Gaussian variable. Since $\BBE[\varepsilon_i Z_j]=0$ for 
all $j\neq i$, and all $(\varepsilon, Z)$ are jointly Gaussian, it 
follows that $\varepsilon_i$ is independent of $Z_{\backslash i}$ as 
claimed. 
\end{proof}

\begin{proof}[Proof of Theorem~\ref{th:CISinDAG}]
Using Lemma~\ref{lem:nonnegprecis}(b) recursively starting with $i=m$ 
we get that $X$ is CIS if and only if 
\[
X_i\;=\;\sum_{j=1}^{i-1} \Lambda_{ij}X_j+\varepsilon_i\quad\mbox{for 
all }i=1,\ldots,m
\]
with $\Lambda_{ij}\geq 0$ and $\varepsilon_i$ independent of 
$X_1,\ldots,X_{i-1}$. This is precisely \eqref{eq:SEM} when applied to 
the full DAG with $j\to i$ for all $j < i$.
\end{proof}

Theorem~\ref{th:CISinDAG} together with Remark~\ref{rem:precfromlam} 
gives the following important algebraic characterization of Gaussian 
CIS distributions. 
\begin{corollary}\label{cor:UU}
The vector $X\sim \cN_m(\mu,\Sigma)$ is CIS if and only if 
$K=UU^\top$ with $U$ upper triangular with positive diagonal and 
non-positive off-diagonal entries.
\end{corollary}

Note that the CIS property relies on the ordering of the variables in 
the vector $X$. The following definition is a natural extention of the 
CIS property; see also \cite{muller2001stochastic}.
\begin{definition}
\label{def:CISrvar2}
If there exists a permutation $\sigma$ of $[m]$ such that  
$(X_{\sigma(1)}, \ldots, X_{\sigma(m)})$ is CIS, then we say $\sigma$ 
is a CIS ordering of $X$. If \emph{for every} permutation $\sigma$ of 
$[m]$ we have that the vector $(X_{\sigma(1)},\ldots,X_{\sigma(m)})$ is 
also CIS, then we say $X$ is conditionally increasing (CI). 
\end{definition}
Interestingly, in the Gaussian case CI equals MTP$_2$ (see Section~\ref{sec:examples}). 
Next, let $G = (V,E)$ be a DAG. A permutation $\sigma$ of $V$ is a 
\textit{topological ordering} if  $a\to b$ implies $\sigma(a) 
< \sigma(b)$. 
It is well-known that if $G$ is a DAG, there exists a permutation of 
$V$ that is a topological ordering. In relation to the structural 
equation model, it is useful to recall that if a DAG is topologically 
ordered then Remark~\ref{rem:precfromlam} takes on a particularly nice 
form with $\Lambda$ lower triangular. Denote by $\CIS_\sigma$ the set 
of all Gaussian distributions such that 
$(X_{\sigma(1)},\ldots,X_{\sigma(m)})$ is CIS. The following result gives an important characterization of Gaussian 
positive DAG dependent distributions $\BBM_+(G)$.
\begin{theorem}
\label{thm:cisandposreg}
For  a DAG $G$ it holds that
\[
\BBM_+(G)\;\;=\;\;\BBM(G)\cap \CIS_\sigma,
\]
where $\sigma$ is \emph{any} topological ordering of $G$.
\end{theorem}
\begin{proof}
We first show $\BBM_+(G) \subseteq \BBM(G) \cap \CIS_\sigma$. The 
inclusion $\BBM_+(G)\subseteq \BBM(G)$ follows by definition. To argue 
that $\BBM_+(G)\subseteq \CIS_\sigma$ let $\widehat G$ be the complete 
DAG whose only topological ordering is $\sigma$. It is clear that 
$\BBM_+(G)\subseteq \BBM_+(\widehat G)$ and  $\BBM_+(\widehat 
G)=\CIS_\sigma$ by Theorem~\ref{th:CISinDAG}. Consequently, 
$\BBM_+(G)\subseteq \BBM(G)\cap \CIS_\sigma$.

To show the opposite inclusion, note that if $X$ has distribution in 
$\BBM(G)$ then the representation \eqref{eq:SEM} holds. Since $X$ is 
$\CIS_\sigma$ and $\sigma$ is a topological ordering, we get from 
Lemma~\ref{lem:nonnegprecis}(b) that the coefficients $\Lambda_{ij}$ 
must be non-negative and so the distribution of $X$ lies in 
$\BBM_+(G)$.
\end{proof}

Although we focus in this paper 
on the Gaussian case, we note that Lemma~\ref{lem:nonnegprecis} 
suggests a more general definition, which is in line with 
\cite{vanderweele2010signed,wellman1990fundamental}. Consider a random 
vector $X$ with values in $\mathcal X=\prod_{i=1}^m \mathcal X_i$ where 
$\mathcal X_i\subseteq \BBR$. We always assume that $X$ admits a density 
function with respect to some product measure on $\mathcal X$.
\begin{definition}
\label{def:generalposdep}
Suppose that $X$ is a distribution that is Markov to a directed acyclic 
graph $G$. Then $X$ is positively DAG dependent with respect to $G$ if, 
for every $i$ the condition survival function
\[
\BBP\Big(\{X_i\geq x_i\}\Big|\bigcap_{j\in \Pa(i)}\{X_j=x_j\}\Big)
\]
is non-decreasing in $(x_j)_{j\in \Pa(i)}$.
\end{definition}
We will use this more general definition to motivate some non-Gaussian 
examples int he following discussion.

\subsection{Motivating examples}  
\begin{figure}[h]
\includegraphics{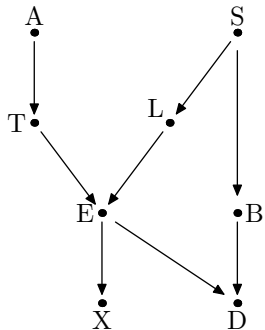}
\caption{The node letters are the parenthetical letters in the list 
above. The variable E represents the logical statement ``Tuberculosis 
(T) or Lung Cancer (L)''.} 
\label{fig:lauritz}
\end{figure}

Positive DAG dependence is often 
present in small, well-designed studies. Some examples of datasets that 
both are well modeled by DAGs and the 
variables in the system are positively correlated, can be found in 
educational research or medical psychology; see, e.g.,  
\cite{afari2013effects,kuipers2019}. There are also two general popular 
datasets, where DAG positive dependence appears naturally. These are 
fictitious datasets that were constructed to mimic real processes. The 
first dataset was introduced in \cite{LS88}. It consists of 
sequences of ``yes'' and ''no'' responses from patients with suspected 
lung disease to the following questions:
\begin{enumerate}
\label{en:list}
\item[(D)] Has shortness-of-breath
\item[(A)] Had a recent trip to Asia
\item[(L)] Has Lung Cancer
\item[(T)] Has Tuberculosis
\item[(E)] Either (T) or (L), or both, are true
\item[(X)] Has a chest X-ray with a positive test
\item[(S)] Is a smoker
\item[(B)] Has Bronchitis
\end{enumerate}
In modeling the relationships of these variables, we take $1$ to be the 
response ``yes'' and $0$ to be ``no'' and use a binary-valued Bayesian 
network illustrated in Figure \ref{fig:lauritz} below to encode the 
relationships between variables, following \cite{LS88}. 

In \cite[Table 1]{LS88}, a ground truth joint distribution was defined 
for this example using the conditional probabilities
\[
\begin{matrix}
\BBP\big( \{\mathrm{A} = 1 \}\big) = .01 & \BBP\big(\{\mathrm{S} = 
1\}\big) = .50 \\
\BBP\big(\{\mathrm{T} = 1\} | \{\mathrm{A} = 1\}\big) = .05 & 
\BBP\big(\{\mathrm{L} = 1\} | \{ \mathrm{S} = 1\}\big) = .10\\
\BBP\big(\{\mathrm{T} = 1\} | \{\mathrm{A} = 0\}\big) = .01 
&\BBP\big(\{\mathrm{L} = 1\} | \{ \mathrm{S} = 0\}\big) = .01 \\
\BBP\big(\{\mathrm{E} = 1\} | \{\mathrm{L} = 1\}\cap \{\mathrm{T} = 
1\}\big) = 1 & \BBP\big(\{\mathrm{B} = 1 \} | \{\mathrm{S} = 1\}\big) = 
.60 \\
\BBP\big(\{\mathrm{E} = 1\} | \{\mathrm{L} = 1\}\cap \{\mathrm{T} = 
0\}\big) = 1 &\BBP\big(\{\mathrm{B} = 1 \} | \{\mathrm{S} = 0\}\big) = 
.30 \\
\BBP\big(\{\mathrm{E} = 1\} | \{\mathrm{L} = 0\}\cap \{\mathrm{T} = 
1\}\big) = 1 & \BBP\big(\{\mathrm{X} = 1 \} | \{\mathrm{E} = 1 \}\big) 
= .98  \\
\BBP\big(\{\mathrm{E} = 1\} | \{\mathrm{L} = 0\}\cap \{\mathrm{T} = 
0\}\big) = 0 & \BBP\big(\{\mathrm{X} = 1 \} | \{\mathrm{E} = 0 \}\big) 
= .05\\
\BBP\big(\{\mathrm{D} = 1\} | \{\mathrm{E} =1\} \cap \{\mathrm{B} = 1\} 
\big) = .90 & \BBP\big(\{\mathrm{D} = 1\} | \{\mathrm{E} = 1\} \cap 
\{\mathrm{B} = 0 \} \big) = .70\\
\BBP\big(\{\mathrm{D} = 1\} | \{\mathrm{E} = 0\} \cap \{\mathrm{B} = 
1\} \big) = .80 & \BBP\big(\{\mathrm{D} = 1\} | \{\mathrm{E} = 0\} \cap 
\{\mathrm{B} = 0\} \big) = .10. 
\end{matrix}
\]
It is clear that the above model is positive dependent with respect to 
the given DAG by inspecting the probabilities directly and checking 
that the condition in Definition~\ref{def:generalposdep} holds.

\begin{figure}[h]
\includegraphics{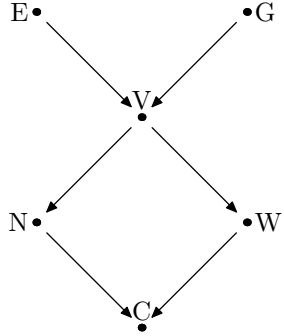}
\caption{The DAG representing the crop dataset from 
\cite{scutari2021bayesian}. The nodes are: E (environmental potential), 
G (genetic potential), V (vegetative organs), N (number of seeds), W 
(seeds mean weight), C (crop).}
\label{fig:crop}
\end{figure}

Another dataset that is used in the context of Gaussian DAGs is the 
crop analysis dataset discussed in Section~2.1 in 
\cite{scutari2021bayesian}. The underlying DAG and the node 
descriptions is given in Figure~\ref{fig:crop}. The dataset assumes the 
following conditional node distributions:
\begin{align*}
    \mathrm{E} &\sim \cN(50, 100)\\
    \mathrm{G} &\sim \cN(50, 100)\\
	\mathrm{V}\;|\;\mathrm{G},\mathrm{E} & \sim \cN(-10.36+0.5 
	\mathrm{G}+0.77\mathrm{E},25)\\
    \mathrm{N}\;|\;\mathrm{V} &\sim \cN(45+0.1\mathrm{V},99)\\
    \mathrm{W}\;|\;\mathrm{V} &\sim \cN(15+0.7\mathrm{V},51) \\
	\mathrm{C}\;|\;\mathrm{N},\mathrm{W} &\sim  \cN(0.3\mathrm{N} + 
	0.3\mathrm{W},39.06)
\end{align*}
Here, again, positive DAG dependence is part of the construction 
because all the conditional means depend positively on the conditioning 
variables. 

\section{Illustrative theoretical examples}
\label{sec:examples}

Denote $\MTP_2$ to be the set of all MTP${}_2$ Gaussians, and $\PA$ to 
be the set of all positively associated Gaussians; see \cite{EPW67} for 
a discussion of association.  In \cite{muller2001stochastic} it is 
shown that for general distributions, the MTP${}_2$ property implies CI 
which in turn implies CIS, and CI is equal to MTP${}_2$ in the Gaussian 
case. Thus, in the Gaussian case, for every permutation $\sigma$ we 
have:
\begin{equation}\label{eq:inclusions}
\MTP_2\quad=\quad\bigcap_\tau \CIS_\tau\quad\subset \quad 
\CIS_\sigma\quad\subset\quad \bigcup_\tau 
\CIS_\tau\quad\subset\quad\PA,    
\end{equation}
where the intersection and the union are taken over all orderings. As 
we will see, even in the Gaussian case, all the inclusions are strict. We first give a simple example which is not MTP${}_2$ but is CIS. 

\begin{example}\label{ex:CISnotMTP}
Consider the upper triangular matrix 
\[
U = \begin{bmatrix}
1 & 0 & -a \\
0 & 1 & -b \\
0 & 0 & 1
\end{bmatrix}
\]
with $a,b>0$. If $K = UU^\top$ is the precision matrix of a Gaussian $X 
= (X_1,X_2,X_3)$, then $X$ is CIS by Corollary~\ref{cor:UU}. However,
\[
K = \begin{bmatrix}
1 & 0 & -a \\
0 & 1 & -b \\
0 & 0 & 1
\end{bmatrix}
\begin{bmatrix}
1 & 0 & 0 \\
0 & 1 & 0 \\
-a & -b & 1
\end{bmatrix}
= \begin{bmatrix}
1 + a^2 & ab & -a \\
ab &  1 + b^2 & -b \\
-a & -b & 1
\end{bmatrix},
\]
which is not an M-matrix, therefore $X$ is not MTP${}_2$. As a 
structural equation model, we may write $X$ as
\begin{align*}
X_1 &= \varepsilon_1\\
X_2 &= \varepsilon_2\\
X_3 &= a X_1 + b X_2 +\varepsilon_3,
\end{align*}
where $(\varepsilon_1,\varepsilon_2,\varepsilon_3)$ is a standard 
$\cN(0,I_3)$ Gaussian. This is a DAG with a v-structure, $1\rightarrow 
3 \leftarrow 2$. Note that $K_{12}>0$ and so $123$ and $213$ are the 
only possible CIS orderings. 
\end{example}

The above example is significant in that it shows that for Gaussian 
distributions, the class of CIS ordered graphical models is 
substantially larger than MTP${}_2$ Gaussians. In particular it is 
known that v-structures cannot occur for MTP${}_2$ graphical models in 
a very general setting \cite{FLSUWZ17}. From this standpoint, it is 
quite appealing to be able to extend from MTP${}_2$ distributions to 
CIS distributions, since v-structures are significant in determining 
Markov equivalence classes, which we discuss in the next section.

Example~\ref{ex:CISnotMTP}  shows that a distribution that is CIS may 
not be CIS with respect to other orderings. In consequence, the 
inclusion $\CIS_\sigma\subset \bigcup_\tau \CIS_\tau$ is also strict 
(unless $m=2$). As a demonstration that the last inclusion in 
\eqref{eq:inclusions} is strict, we give the following example which is 
a positive associated Gaussian where \textit{no reordering} of $X$ is 
CIS.

\begin{example}
Let $X$ be a centered Gaussian with covariance
\[
\Sigma = \begin{bmatrix}
5 & 4 & 7 & 8 \\
4 & 9 & 8 & 7 \\
7 & 8 & 11 & 11 \\
8 & 7 & 11 & 14
\end{bmatrix}.
\]
Since all entries of $\Sigma$ are positive, by 
\cite{pitt1982positively}, $X$ is a positive associated Gaussian. 
However,
\[
K = 
\begin{bmatrix}
94 & 25 & -55 & -23 \\
25 &  7 & -15 &  -6 \\
-55 & -15 & 33 & 13 \\
-23 & -6 & 13 & 6
\end{bmatrix},
\]
since each row of the above matrix  has a positive off-diagonal entry 
it follows that there is no $j \in [4]$ such that $\BBE[X_j | 
X_{\backslash j}]$ is a non-decreasing function in $X_{\backslash j}$, 
from which we conclude that there is no CIS ordering of $X$.
\end{example}

The next result studies the relation between $\CIS_\sigma$ models.
\begin{proposition}
\label{prop:ciseq}
Suppose $X=(X_1,\ldots,X_m)$ has a Gaussian distribution. If $m=2$ then 
$(X_1,X_2)$ is {\rm CIS} if and only if $(X_2,X_1)$ is {\rm CIS}. If 
$m\geq 3$  then ${\rm CIS}_\sigma={\rm CIS}_{\sigma'}$ if and only if 
$\sigma(k)=\sigma'(k)$ for $k=3,\ldots,m$. 
\end{proposition}
\begin{proof}
The bivariate case follows because $(X_1,X_2)$ is CIS if and only if 
$\Cov(X_1,X_2)\geq 0$, which is symmetric in $(X_1,X_2)$. Suppose 
$m\geq 3$. The ``if'' implication follows directly from the definition 
and  from the $m=2$ case. For the ''only if`` implication assume with 
no loss in generality that $\sigma'={\rm id}$. We construct a 
distribution that lies in $\CIS_{\rm id}$ and show that it lies in 
$\CIS_{\sigma}$ if and only if $\sigma={\rm id}$ or 
$\sigma=(2,1,3,\ldots,m)$. Let $U$ be an upper triangular matrix of the 
form
\[
U\;=\;\begin{bmatrix}
1 & 0 & -1 & -2 & \cdots & -(m-3)&-(m-2)\\
0 & 1 & -1 & -1 & \cdots & -1&-1\\
0 & 0 & 1 & -1 & \cdots & -1&-1\\
\vdots & \vdots &  & \ddots & \ddots &  & \\
0 & 0 & 0 & 0 & \cdots & 1&-1\\
0 & 0 & 0 & 0 & \cdots & 0& 1
\end{bmatrix}.
\]
The distribution we construct has covariance $\Sigma$ such that 
$K=\Sigma^{-1}=U U^\top$. Since all the upper off-diagonal entries are 
non-positive, this distribution is $\CIS_{\rm id}$. Denote the rows of 
$U$ by $U_1,\ldots,U_m$. Note that 
\[
U_1^\top U_2\;>\;\cdots \;>\;U_1^\top U_{m-1}\;=\;1\;>\;0,
\]
and so $K_{1i}>0$ for all $i=2,\ldots,m-1$. This shows that every CIS 
ordering of this random vector must have $m$ as the last index. If $m= 
3$, then we are done. If $m\geq 4$, consider the marginal distribution 
over $A=\{1,\ldots,m-1\}$. Because $U$ is upper triangular, we get that 
$(\Sigma_{A,A})^{-1}=U_{A,A}U_{A,A}^\top$. Note that $U_{A,A}$ has the 
same form as $U$ but with $m-1$ replacing $m$. Thus, by the same 
argument as above,
$$
(\Sigma_{A,A})^{-1}_{1i}>0\quad\mbox{for all }i=2,\ldots,m-2.
$$
This shows that every CIS ordering of our constructed distribution must 
have $m-1$ as the penultimate index. If $m=4$, we are done. If $m\geq 
5$, take $A\setminus \{m-1\}$ as the new $A$ and proceed as above. In 
this way we show that for this distribution $\sigma$ is a CIS ordering 
only if $\sigma(k)=k$ for $k=3,\ldots,m$.
\end{proof}

There are qualitative properties of CIS distributions that contrast with 
MTP${}_2$ distributions. It is known (\cite[Proposition 3.2]{KR80}) 
that if $X$ is MTP${}_2$ distributed then any marginal distribution of 
$X$ also satisfies the MTP${}_2$ property. The next example shows that 
a Gaussian CIS random variable does not satisfy this property.

\begin{example}
Let $X=(X_1,X_2,X_3,X_4)$ be a centered Gaussian with covariance
\[
\Sigma=\begin{bmatrix}
\tfrac14 & \tfrac14 & \tfrac34 & \tfrac{29}{16}\\[.1cm]
\tfrac{1}{4} & \tfrac{5}{4}& \tfrac{7}{4}& \tfrac{77}{16}\\[.1cm]
\tfrac{3}{4}& \tfrac{7}{4}& \tfrac{17}{4}& \tfrac{167}{16}\\[.1cm]
\tfrac{29}{16}& \tfrac{77}{16}& \tfrac{167}{16}& \tfrac{1737}{64}
\end{bmatrix}.
\]
It can be checked directly that $(X_1,X_2,X_3,X_4)$ is CIS. However, 
the inverse of $\Sigma_{134}$ is
\[
\begin{bmatrix}
\tfrac{205}{24} & -\tfrac{23}{12} &  \tfrac{1}{6}\\[.1cm]
-\tfrac{23}{12} &  \tfrac{14}{3} &  -\tfrac{5}{3}\\[.1cm]
\tfrac{1}{6} &  -\tfrac{5}{3} &  \tfrac{2}{3}
\end{bmatrix}.
\]
Since the last row of this matrix has a positive off-diagonal entry, we 
conclude that $(X_1,X_3,X_4)$ is not CIS. 
\end{example}

However, the following result, which follows immediately from the definition, shows that certain conditionals and marginals preserve the CIS property. 
result . 
\begin{proposition}
Let $X$ be a CIS distributed centered Gaussian. Then the following 
distributional properties hold:
\begin{enumerate}
\item The conditional distribution of $(X_{k+1},\ldots, X_m)$ given
$(X_1,\ldots, X_k)$ is CIS.
\item The vector $(X_1,\ldots,X_k)$ is CIS for every $1 \leq k \leq m$.
\end{enumerate}
\end{proposition}
 
Theorem~\ref{thm:cisandposreg} shows a relation between CIS orderings 
and positive DAG dependence. The following example describes a complication that can arise. Namely, we consider a DAG 
whose possible topological orderings are $1(23)4$, the union of all 
$\CIS$ orderings for all Markov equivalent DAGs is $(123)4$, but 
\emph{for special} distributions in the model it is possible that 
$4321$ is a valid $\CIS$ ordering.

\begin{example}
Consider the DAG model defined by the upper triangular matrix 
\[
U = \begin{bmatrix}
1 & -a & -b &0 \\
0 & 1 & 0 & -c\\
0 & 0 & 1 & -d\\
0 & 0 & 0 &1
\end{bmatrix}
\]
with $a,b,c,d>0$. The Markov equivalence class of this DAG consists of the following three DAGs.
\begin{figure}[h]
\includegraphics{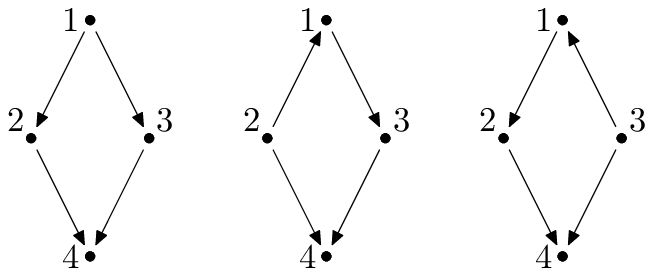}
\end{figure}

\noindent The corresponding precision matrix is given by
\[
K=UU^\top=\begin{bmatrix}
1+a^2+b^2 &	-a	&-b	&0\\
-a	&1+c^2&	cd	&-c\\
-b	&cd	&1+d^2&	-d\\
0	&-c	&-d&	1
\end{bmatrix}
\]
Since $K_{23}=cd>0$, it is clear that any CIS ordering has $1$ or $4$ 
as the last element and $1$ is actually possible. Since $(X_1,X_2,X_3)$ 
is always CI, we conclude that all orderings $(123)4$ are 
CIS\footnote{The notation $(123)$ stands for any permutation of these 
three.}. By direct computation we see that for $\Sigma=K^{-1}$,
\[
(\Sigma_{234})^{-1}\;=\;
\begin{bmatrix}
\frac{(a^2+1)c^2+b^2(c^2+1)+1}{1+a^2+b^2} & 
\frac{a^2cd-ab+(b^2+1)cd}{1+a^2+b^2} &-c\\
\frac{a^2cd-ab+(b^2+1)cd}{1+a^2+b^2} & 
\frac{a^2(d^2+1)+(b^2+1)d^2+1}{1+a^2+b^2} & -d\\ -c & -d & 1
\end{bmatrix}.
\]
In particular, if $a,b$ are sufficiently large and $c,d$ are 
sufficiently small such that $a^2cd-ab+(b^2+1)cd\leq 0$, we also have 
that $(X_2,X_3,X_4)$ is CI. In this case, each ordering $(234)1$ is 
also a CIS ordering. 
Thus the CIS orderings are of the form $1(23)4$, $4(23)1$ and 
$(23)(14)$.
Note that only the DAG with topological ordering $1(23)4$ is in the 
Markov equivalence class, while the DAGs with topological ordering $4(23)1$ and (23)(14) are
not. 
This shows 
that the set of all CIS orderings contains only limited information 
about the underlying DAG. 
\end{example}

The situation is not always that complicated. In 
Proposition~\ref{prop:belief-net} we will show that there is a large 
class of interesting hierarchical networks for which the possible 
$\CIS$ orderings are exactly the topological orderings.  

Another property worth noting is that the space of MTP${}_2$ Gaussian 
distributions amounts to the M-matrix constraint on $K$, which is 
convex in $K$. We can show that the space of $K$ for which $X$ is a CIS 
Gaussian is not convex in $K$.

\begin{example}
Let $K_1$ and $K_2$ be the precision matrices
\[
\begin{split}
K_1 &= 
\begin{bmatrix}
 1 & -1 & -1 & -4 \\
 0 & 1 & 0 &  0 \\
 0 & 0 &  1 & -3 \\
 0 & 0 &  0 &  1
 \end{bmatrix}
 \begin{bmatrix}
 1 & -1 & -1 & -4 \\
 0 & 1 & 0 &  0 \\
 0 & 0 &  1 & -3 \\
 0 & 0 &  0 &  1
 \end{bmatrix}^\top 
=
\begin{bmatrix}
19 & -1 & 11 & -4\\
-1 & 1 & 0 & 0\\
11 & 0 & 10 & -3\\
-4 & 0 & -3 & 1
\end{bmatrix},
\\
K_2 &=
\begin{bmatrix}
1 & -1 & 0 & 0\\
0 & 1 & -1 & 0\\
0 & 0 & 1 & -1\\
0 & 0 & 0 & 1
\end{bmatrix}
\begin{bmatrix}
1 & -1 & 0 & 0\\
0 & 1 & -1 & 0\\
0 & 0 & 1 & -1\\
0 & 0 & 0 & 1
\end{bmatrix}^\top
=	
\begin{bmatrix}
2 & -1 & 0 & 0 \\
-1 & 2 & -1 & 0 \\
0 & -1 & 2 & -1\\
0 & 0 & -1 & 1
\end{bmatrix}.
\end{split}
\]
Clearly, by Corollary~\ref{cor:UU}, we must have that Gaussians with
precision matrices $K_1$ or $K_2$ are CIS ordered. However consider the
sum
\[
K = K_1 + K_2 = 
\begin{bmatrix}
21 & -2 & 11 &-4 \\
-2 & 3 & -1 & 0 \\
11 & -1 & 12 & -4 \\
-4 & 0 & -4 & 2
\end{bmatrix}.
\]
Then if $\Sigma = K^{-1}$, by the Schur Complement formula
\[
\big(\Sigma_{[3],[3]}\big)^{-1} = K_{[3],[3]} - \frac{K_{[3],4} 
K_{4,[3]}}{K_{4,4}}
=
\begin{bmatrix}
13 & -2 & 3\\
-2 &  3 & -1\\
3 & -1 & 4
\end{bmatrix}
\]
which means that if $X$ has covariance $\Sigma$, then 
$\BBE[X_3|X_1,X_2]$ is not a non-decreasing function of $X_1$ and $X_2$ 
due to the third row of the above matrix having an off-diagonal that is 
positive --- the same will be true if we were to replace $K$ by 
$\frac{K}{2}$, which implies that the convex combination of $K_1$ and 
$K_2$ does not stay in the class of precision matrices of CIS 
Gaussians.
\end{example}

The above example shows that even if we assume that a Gaussian is CIS 
ordered under a known permutation $\sigma$, we do not have convexity in 
the space of precision matrices $K$ that parameterize this model. In 
Section~\ref{sec:choleskestimat}, we show that there is a broad class 
of models, under which Gaussians that are CIS for a known permutation 
$\sigma$ are included, for which computing the MLE is a convex 
optimization problem. While the results of 
Section~\ref{sec:choleskestimat} may be familiar to many practitioners, 
we did not find a direct reference and thought it worthwhile to specify 
these models explicitly. Most importantly for us however, is that 
computationally, once a CIS ordering is known, calculating the MLE for 
a CIS Gaussian can be done with similar efficiency as restricting to 
the MTP${}_2$ class. 
 
\section{Markov equivalence for CIS models}
\label{sec:markovandcis}

One of the most fundamental limitations of Bayesian networks is that 
two different DAGs may represent the same conditional independence 
model, in which case we say that they are Markov equivalent. We recall 
the following classical result \cite{verma2022equivalence} that uses 
the concept of a skeleton, which for a DAG $G$ is the undirected graph 
obtained from $G$ by forgetting the directions of all arrows.
\begin{theorem}
\label{thm:markovequiv}
Two DAGs $G$ and $H$ are Markov equivalent if and only if they have the 
same skeleton and v-structures. For a Gaussian $X$, we have 
$\BBM(G)=\BBM(H)$ if and only if $G$ and $H$ are Markov equivalent.
\end{theorem}
If $G$ is a DAG then by $[G]$ we denote the set of all DAGs Markov 
equivalent to $G$. There is another useful characterization of Markov 
equivalence proposed in~\cite{chickering1995transformational}, which describes elementary 
operations on DAGs that transform a DAG into a Markov equivalent DAG in 
such a way that $G$ can be transformed to any graph in $[G]$ by a 
sequence of these elementary operations. This elementary operation is 
given by flipping the arrow $i\to j$ whenever the pair $i,j$ satisfies 
$\{i\}\cup \Pa(i)=\Pa(j)$. More specifically, given a DAG $G$ over $V$, 
we say that an arrow $i\to j$ in $G$ is \textit{covered} if the graph 
$H$ obtained from $G$ by replacing $i\to j$ with $j\to i$ is also 
acyclic, and also ${\rm Pa}(i)={\rm Pa}(j)\setminus \{i\}$. 
result of \cite{chickering1995transformational} states:
\begin{theorem}\label{th:chickering}
We have $H\in [G]$ if and only if $H$ can be obtained from $G$ by a 
sequence of flips of covered arrows. 
\end{theorem}
 
We say that $H$ is CIS-Markov equivalent to $G$ if 
$\BBM_+(G)=\BBM_+(H)$. We offer a similar characterization of 
CIS-Markov equivalence. An edge $i\to j$ is \emph{trivially covered} if  
${\rm Pa}(i)\cup \{i\}={\rm Pa}(j)=\{i\}$.
\begin{theorem}\label{th:carolineconj}
For a Gaussian $X$, we have $\BBM_+(G)=\BBM_+(H)$ if and only if $H$ 
can be obtained from $G$ by a sequence of flips of trivially covered 
arrows. 
\end{theorem}
Note that when $G$ is a complete DAG (with all possible $\binom{m}{2}$ 
edges), then $\BBM_+(G)=\CIS_\sigma$, where $\sigma$ is the unique 
topological ordering of this DAG. This shows that 
Theorem~\ref{th:carolineconj} generalizes Proposition~\ref{prop:ciseq}. 
\begin{proof}
For the ``if'' part, it is enough to consider the case when $H$ is 
obtained from $G$ by a single flip of a trivially covered pair $i\to 
j$. By Theorem~\ref{th:chickering}, $\BBM(G)=\BBM(H)$. Since $i$ has no 
parents and it is the only parent of $j$, there is a permutation 
$\sigma$ with $\sigma(1)=i$, $\sigma(2)=j$ that forms a topological 
ordering of $G$. Moreover,  the permutation $\sigma'$ obtained from 
$\sigma$ by swapping $i$ and $j$ is a topological ordering of $H$. By 
Proposition~\ref{prop:ciseq}, $\CIS_\sigma=\CIS_{\sigma'}$. By 
Theorem~\ref{thm:cisandposreg}, $\BBM_+(G)=\BBM_+(H)$.

To show the ``only if'' part, first note that if $\BBM_+(G)=\BBM_+(H)$ 
then necessarily $\BBM(G)=\BBM(H)$. Via a contrapositive 
argument, suppose $G$ and $H$ are Markov equivalent but $H$ is not 
obtained from $G$ by a sequence of trivially covered edge flips. 
This means that there exists an arrow $i\to j$ in $G$ and $k$ with 
$k\in \Pa_G(i)\cap \Pa_G(j)$ such that $i\leftarrow j$ in 
$H$. To get a contradiction, it is enough to construct a distribution 
in $\BBM_+(G)$ such that in every CIS ordering $j$ must come after $i$.

Let $\sigma$ be a topological ordering of $G$. Without loss of 
generality assume $\sigma={\rm id}$ and let $i,j,k$ be as above. In 
particular, $1\leq k<i<j\leq m$. Let $U$ be upper triangular such that 
$U_{ll}=1$ for all $l=1,\ldots,m$, $U_{ij}=-1$, $U_{kj}=-1$ and $U$ is 
zero otherwise. Note that by the above, this $U$ corresponds to a 
distribution in  $\BBM_+(G)$ where some of the edges in $G$ have zero 
coefficients. We will show that for any $A$ containing $\{i,j,k\}$, 
neither $i$ nor $k$ can be the last one in a CIS ordering. To show 
this, note that $U_{A,A^c}=0$, $U_{A^c,A}=0$, and $U_{A^c,A^c}=I$. It 
follows that 
$$
(\Sigma_{A,A})^{-1}\;=\;U_{A,A}U_{A,A}^\top
$$
and so $(\Sigma_{A,A})^{-1}_{ik}=1>0$ showing that neither $i$ nor $k$ 
can be the last element in any CIS ordering of $X_A$. Using this 
recursively, starting from $A=\{1,\ldots,m\}$, we conclude that $j$ 
must appear after $i,k$ in every CIS ordering.
\end{proof}

In Gaussian Bayesian networks the crucial observation is that if the 
Markov equivalence classes $[G]$ and $[H]$ are not equal then the 
Gaussian models $\BBM(G)$ and  $\BBM(H)$ intersect at a measure zero 
set (we think about the geometry of these models as embedded in the 
space of covariance matrices). This means that for almost all ground-truth models
we can learn the equivalence classes from the data. The analogous 
statement is unfortunately not true for CIS-Markov equivalence 
classes. For example, if $m=3$, the following two graphs lie in the 
same Markov equivalence class and different CIS-Markov equivalence 
classes
\begin{figure}[h]
\includegraphics{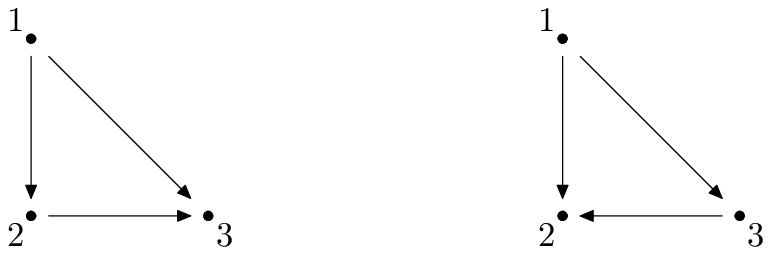}
\end{figure}

\noindent The intersection of $\BBM_+(G)$ and $\BBM_+(H)$ for these two graphs 
has full dimension and it contains the set of all inverse $M$-matrices.

\begin{lemma}
Suppose the distribution of $X$ lies in $\BBM_+(G)$. Suppose that there 
exists $k$ such that  $i\to k\leftarrow j$ is a v-structure in $G$ and 
suppose that $K_{ij}\neq 0$ (this holds generically).  Then 
no CIS ordering of $X$ finishes with $i$ or $j$.
\end{lemma}
\begin{proof}
Without loss of generality assume that the trivial ordering 
$1<2<\ldots<m$ is a topological ordering of $G$. In this case the 
matrix $\Lambda$ in \eqref{eq:SEM} is lower triangular. Then let 
$K=UU^\top$, with $U$ upper triangular, be the precision matrix of 
$X$. By Remark~\ref{rem:precfromlam} we have $U=(I-\Lambda)^\top 
D^{-1/2}$ and so for $i\neq j$ have $U_{uv}\leq 0$ if $u\to v$ in $G$ 
and $U_{uv}=0$ otherwise. We have
$$
K_{ij}\;=\;\sum_l U_{il}U_{jl}\;=\;\sum_{l\in \Ch(i)\cap \Ch(j)} 
U_{il}U_{jl}.
$$
This expresses $K_{ij}$ as a sum of non-negative terms. Since this sum 
is non-zero by assumption, it must be strictly positive and so, neither 
$i$ nor $j$ can be the last ones in a CIS ordering.
\end{proof}
As a corollary we get the following result.
\begin{proposition}\label{prop:belief-net}
Consider a DAG $G$ consising of $k$ layers $V_1,\ldots,V_k$ such 
that:
\begin{enumerate}
\item $V=V_1\sqcup\cdots\sqcup V_k$,
\item only arrows from $V_{i}$ to $V_{i+1}$ are allowed in $G$,
\item $|V_i|\geq 2$ for all $i=1,\ldots,k-1$ (only the last layer may 
contain one node),
\item every $v\in V_i$ for $i=1,\ldots,k-1$ is contained in a 
v-structure (as a parent).
\end{enumerate}
If the distribution of $X$ lies in $\BBM_+(G)$ and $K_{ij}\neq 0$ 
unless $i,j\in V_k$ (this holds generically), then the only 
possible CIS orderings of $X$ are $(V_1)\cdots(V_k)$, where the 
notation $(V_i)$ means that the vertices in $V_i$ can be ordered in an 
arbitrary way. In particular, any possible CIS ordering of $X$ is a 
topological ordering of $G$.
\end{proposition}

\section{Maximum likelihood estimation in $\BBM_+(G)$}
\label{sec:choleskestimat}

In this section we show that maximum likelihood estimation in the model 
$\BBM_+(G)$ for a given $G$ is straightforward and amounts to solving a 
convex optimization problem. 
Consider a Gaussian vector $X\sim \cN_m( 0,\Sigma)$ and let 
$K=\Sigma^{-1}$. Since $K$ is positive definite, by \cite[Corollary 
3.5.6]{HJ13} we have that there exists a unique upper triangular matrix 
$U$ whose diagonals are all 1, and a diagonal matrix $D$ with strictly 
positive diagonals such that $K = UDU^\top$. Moreover, the relation 
between $K$ and the pair $(D,U)$ is one-to-one. Equivalently, we obtain 
the stochastic representation 
\begin{equation}\label{eq:SR}
X=\Lambda X+\ve,
\end{equation}
where $\Lambda=(I_m-U)^\top$ is lower triangular with zero diagonals, 
and $\ve\sim \cN_m( 0,D^{-1})$.

\begin{definition}
Let $\mathcal L_i\subseteq \R^{i}$ be sets for each $i=1,\ldots,m-1$. A 
\emph{Cholesky factor model} consists of all Gaussian distributions 
such that the inverse-covariance matrix satisfies $K=UDU^\top$ with  
$D$ a diagonal matrix and $U=(I_m-\Lambda)^\top$ with 
\[
\Lambda_i:=(\Lambda_{i,1},\ldots,\Lambda_{i,i-1})\in \mathcal 
L_{i-1}\quad \mbox{for}\quad i=2,\ldots,m.
\]
\end{definition}
\begin{remark}
In the case that $\cL_i = [0,\infty)^{i}$, we recover the CIS model on 
$X$. If $\cL_i = \BBR^{i}$ we simply have the space of all covariance 
matrices.
\end{remark}

\begin{remark}
If the DAG $G$ is known, we can always assume without loss of 
generality that the ${\rm id}$ permutation is a topological ordering of 
$G$. In other words, the matrix $\Lambda$ in 
Remark~\ref{rem:precfromlam} is lower triangular. Thus $\BBM(G)$ is a 
Cholesky factor model with the support of $\Lambda_i$ equal to the 
parent set ${\rm Pa}(i)$. The model $\BBM_+(G)$ is obtained by 
additional non-negativity constraints.
\end{remark}

If we want to make the constraints on $\Lambda$ explicit we denote the 
model by $F(\mathcal L_1,\ldots,\mathcal L_{m-1})$. Maximum likelihood 
estimation for such models links to the problem of least squares estimation in linear regression as follows.  
Given $n$ independent observations of $X$ from this model, we stack 
them in the matrix $\BX\in \R^{n\times m}$. We denote by 
$\mathbf{x}_1,\ldots,\mathbf{x}_m$ the columns of 
$\mathbf{X}$ and by $\mathbf{Z}_i:=\BX_{[n],[i-1]}$ the 
$\R^{n\times (i-1)}$ matrix obtained from the first $i-1$ columns of 
$\BX$. 

\begin{theorem}\label{th:cfmle}
If $(\hat D,\hat \Lambda)$ is the maximum likelihood estimator for a Cholesky factor model $F(\mathcal L_1,\ldots, 
\mathcal L_{m-1})$, then each $\hat \Lambda_i$ for $i=2,\ldots,m-1$ is 
given as a minimizer of the quadratic problem
$$
{\rm minimize}\;\;\tfrac1n \|\mathbf{x}_i-\mathbf{Z}_i\Lambda_i^\top\|^2\qquad\mbox{subject to }\Lambda_i\in\mathcal L_{i-1}\subseteq \R^{i-1}. 
$$
Moreover, 
$$
\hat D_{ii}\;=\;n\|\mathbf{x}_i-\mathbf{Z}_i\hat \Lambda_i^\top\|^{-2}
$$
for all $i=1,\ldots,m$.
\end{theorem}
\begin{proof}
We have $K=(I_m-\Lambda)^\top D (I_m-\Lambda)$, where $\Lambda$ is 
strictly lower triangular with $\Lambda_i\in \mathcal L_{i-1}$ for 
$i=2,\ldots,m$. As before, set $U=(I_m-\Lambda)^\top$. Since 
$\det(U)=1$ and $D$ is diagonal, the corresponding log-likelihood 
function $\log\det(K)-\tfrac{1}{n}{\rm tr}(\mathbf{X}^\top\mathbf{X} K 
)$ can be written as 
\begin{equation}\label{eq:llike}
\sum_{i=1}^m \log D_{ii}-\frac1n\sum_{i=1}^m {D_{ii}}((\BX U)^\top \BX 
U)_{ii}.
\end{equation}
The expression $((\BX U)^\top \BX U)_{ii}$ is simply the squared-norm 
of the $i$-th column of $\BX U$, which is equal 
to 
$$
\mathbf{x}_i-\sum_{j=1}^{i-1}\Lambda_{ij}\mathbf{x}_j\;=\;\mathbf{x}_i 
- \mathbf{Z}_i \Lambda_i^\top.
$$
Thus, maximizing \eqref{eq:llike} is equivalent to minimizing
\begin{equation}\label{eq:llike2}
-\sum_{i=1}^m \log D_{ii}+\sum_{i=1}^m \frac{D_{ii}}{n}\|\mathbf{x}_i- 
\mathbf{Z}_i\Lambda_i^\top\|^2.
\end{equation}
The $i$-th squared term in \eqref{eq:llike2} depends only on 
$\Lambda_i$. This means that minimizing \eqref{eq:llike2} in a 
Cholesky factor model can be done term by term. Once 
the optimizer for $\Lambda$ is found, $D$ can be handled in a 
straightforward way. 
\end{proof}

Theorem~\ref{th:cfmle} gives also a simple condition on the existence 
of the MLE. 
\begin{proposition}\label{prop:MLEexistence}
The MLE in Theorem~\ref{th:cfmle} exists if and only if each set 
$\cL_i$ is closed and for every $i=1,\ldots,m-1$, 
$$
\mathbf{x}_i\;\notin \;\{\mathbf{Z}_i\Lambda_i^\top:\;\Lambda_i\in 
\cL_i\}.
$$
In particular, if there are subsets $A_i \subseteq [i-1]$ such that  
$\cL_i={\rm span}\{\mathbf{x}_j: j\in A_i\}$, then the MLE exists with 
probability 1 as long as $n\geq \max_i |A_i|$. 
\end{proposition}
It is now straightforward to compute the optimal value for the 
log-likelihood. 
\begin{corollary}
If the MLE exists, then the optimal value of the log-likelihood is 
$$
-\sum_{i=1}^m \log\left(\frac1n \|\mathbf{x}_i-\mathbf{Z}_i\hat 
\Lambda_i^\top\|^2\right)-m.
$$
\end{corollary}

Recall that in the linear regression problem, given a vector 
$\mathbf{x}_i\in \R^n$ and the matrix $\mathbf{Z}_i\in \R^{n\times 
(i-1)}$, the least squares estimator is given precisely as the 
minimizer of $\|\mathbf{x}_i-\mathbf{Z}_i\theta\|^2$ over $\theta\in 
\R^{i-1}$. If this minimizer is unique, it is given by the well-known 
formula 
\begin{equation}\label{eq:ols}
\hat\theta=(\mathbf{Z}_i^\top 
\mathbf{Z}_i)^{-1}\mathbf{Z}_i^\top\mathbf{x}_i.
\end{equation}
If $\mathbf{Z}_i$ does not have full column rank, the optimum is 
obtained over an affine space. Replacing the inverse above with the 
pseudo-inverse gives the solution with the smallest norm. The following 
result follows almost immediately.
\begin{proposition}
If the constraints $\mathcal L_1,\ldots, \mathcal L_{m-1}$ are all 
linear, then the MLE $(\hat D,\hat \Lambda)$ in the Cholesky factor 
model $F(\mathcal L_1,\ldots, \mathcal L_{m-1})$ can be given in closed 
form. 
\end{proposition}

\section{Finding a CIS ordering}
\label{sec:cisordrecov}

Having established that the MLE can be easily computed in $\BBM_+(G)$ 
for any fixed $G$, we now explore the harder problem of estimating 
$\Sigma$ knowing that the distribution of $\mathbf X$ lies in  
$\BBM_+(G)$ for \emph{some} $G$. {By Theorem~\ref{thm:cisandposreg}, 
$\BBM_+(G)=\BBM(G)\cap \CIS_\sigma$ for any topological ordering 
$\sigma$ of $G$. Thus, if we know a topological ordering of $G$ the 
problem can be solved by running regressions in 
the order given by $\sigma$ and adding a LASSO penalty to learn a sparse 
representation.} 
Using the fact that $\BBM_+(G)=\BBM(G)\cap 
\CIS_\sigma$ for any topological ordering $\sigma$ of $G$, we do not need to search over all orderings but can 
restrict ourselves to CIS orderings for the underlying distribution. In 
this section, we show that these can be efficiently recovered.

\subsection{Recovering a CIS ordering in the population case}
In the following, we provide an algorithm that, given $K$, recovers a CIS ordering given that such an ordering exists. The algorithm is based on the
following lemma.

\begin{lemma}
\label{thm:reorderingaCIS}
Suppose $X$  is a CIS $m$-variate Gaussian. Suppose there exists $k\in 
[m-1]$ such that $K_{k,\backslash k}\leq 0$. Then 
$(X_{1},\ldots,X_{k-1},X_{k+1},\ldots,X_{m},X_{k})$ is CIS.
\end{lemma}
\begin{proof}
Recalling Lemma~\ref{lem:nonnegprecis}, we have that $X$ being a 
centered CIS ordered Gaussian is equivalent to 
\[
\BBE[X_j | X_{[j-1]}] = - 
\frac{\big(\Sigma_{[j],[j]}\big)^{-1}_{j,[j-1]}}{\big(\Sigma_{[j],[j]}\big)^{-1}_{j,j}} 
X_{[j-1]}
\]
being a non-decreasing function in $(X_1,\ldots, X_{j-1})$ for all $j 
\in [m]$. We only need to check that the functions
\begin{align}
\label{eqn:middle}
\BBE&[X_j | X_{[j-1]\backslash\{k\}}] \qquad j=k+1,\ldots, m,\\
\label{eqn:end}
\BBE&[X_k | X_{\backslash k}]
\end{align}
are non-decreasing in their arguments, which, at least for the second 
function, follows automatically by the assumption $K_{k,\backslash 
k}\leq 0$. We now proceed by an induction argument starting from $j=m$ 
working downward, to prove that the functions~\eqref{eqn:middle} are 
all non-decreasing in their arguments. We have
\[
\BBE\big[X_m | X_{[m-1]\backslash \{k\}}\big] = 
-\frac{\big(\Sigma_{\backslash k,\backslash 
k}\big)^{-1}_{m,[m-1]\backslash \{k\}}}{\big( \Sigma_{\backslash k, 
\backslash k}\big)^{-1}_{m,m}} X_{[m-1]\backslash{\{k\}}},
\]
then the Schur complement formula gives the following two statements:
\[
\begin{split}
\big(\Sigma_{\backslash k,\backslash k}\big)^{-1} &= 
K_{\backslash k,\backslash k} - \frac{K_{\backslash k, k} 
K_{k,\backslash k}}{K_{k,k}},\\
\big(\Sigma_{[m-1],[m-1]}\big)^{-1} &= 
K_{[m-1],[m-1]} - \frac{K_{[m-1], m} K_{m,[m-1]}}{K_{m,m}}.
\end{split}
\]
By our assumption $K_{k,\backslash k} \leq 0$, we have that 
$\frac{K_{\backslash k,k}K_{k,\backslash k}}{K_{k,k}}$ is a 
non-negative rank-one matrix. Similarly, since $X_m$ is the last in a 
CIS ordering of $X$, we have $K_{m,\backslash m} \leq 0$ and 
$\frac{K_{[m-1],m}K_{m,[m-1]}}{K_{m,m}}$ is a non-negative matrix. It 
follows then that
\begin{equation}
\label{eqn:inductionbase}
\begin{split}
\big(\Sigma_{\backslash k,\backslash k}\big)^{-1}_{m,[m-1]\backslash \{k\}} &\leq 0,\\
\big(\Sigma_{[m-1],[m-1]}\big)^{-1}_{k,[m-1]\backslash \{k\}} &\leq 0.
\end{split}
\end{equation}
The first inequality in~\eqref{eqn:inductionbase} implies that the 
function in equation~\eqref{eqn:middle} for $j=m$ is non-decreasing in 
its arguments.
Our induction hypothesis is that for some $j^*\geq k+1$ we have shown
for every $j=j^*+1,\ldots, m$, that 
\begin{equation}
\label{eqn:inductionequation}
\begin{split}
\big(\Sigma_{[j]\backslash \{k\}, [j]\backslash \{k\}}\big)^{-1}_{j,[j-1]\backslash \{k\}} &\leq 0,\\
\big(\Sigma_{[j-1], [j-1]}\big)^{-1}_{k,[j-1]\backslash \{k\}} &\leq 0.
\end{split}
\end{equation}
We will now prove that both of these inequalities are true for 
$j= j^*$ as well. By the second inequality 
in~\eqref{eqn:inductionequation} (setting $j = j^*+1$), and the fact 
that $X$ is CIS ordered, we have that
\begin{equation}
\label{eqn:inductioningred}
\begin{split}
\big(\Sigma_{[j^*],[j^*]}\big)^{-1}_{k,[j^*]\backslash\{k\}} &\leq 0,\\
\big(\Sigma_{[j^*],[j^*]}\big)^{-1}_{j^*,[j^*-1]} &\leq 0.
\end{split}
\end{equation}
The Schur complement formula implies the following two equalities
\begin{multline*}
\big(\Sigma_{[j^*]\backslash\{k\},[j^*]\backslash\{k\}}\big)^{-1} =\\
\big(\Sigma_{[j^*],[j^*]}\big)^{-1}_{[j^*]\backslash\{k\},[j^*]\backslash\{k\}}
- \frac{\big(\Sigma_{[j^*],[j^*]}\big)^{-1}_{[j^*]\backslash\{k\},k}\big(\Sigma_{[j^*],[j^*]}\big)^{-1}_{k,[j^*]\backslash\{k\}}}{\big(\Sigma_{[j^*],[j^*]}\big)^{-1}_{k,k}},
\end{multline*}
and
\begin{multline*}
\big(\Sigma_{[j^*-1],[j^*-1]}\big)^{-1} =\\
\big(\Sigma_{[j^*],[j^*]}\big)^{-1}_{[j^*-1],[j^*-1]}
- \frac{\big(\Sigma_{[j^*],[j^*]}\big)^{-1}_{[j^*-1],j^*}\big(\Sigma_{[j^*],[j^*]}\big)^{-1}_{j^*,[j^*-1]}}{\big(\Sigma_{[j^*],[j^*]}\big)^{-1}_{j^*,j^*}}.
\end{multline*}
By the inequality of equation~\eqref{eqn:inductioningred}, it follows
that 
\[
\begin{split}
\frac{\big(\Sigma_{[j^*],[j^*]}\big)^{-1}_{[j^*]\backslash\{k\},k}\big(\Sigma_{[j^*],[j^*]}\big)^{-1}_{k,[j^*]\backslash\{k\}}}{\big(\Sigma_{[j^*],[j^*]}\big)^{-1}_{k,k}}&\geq 0,
\\
\frac{\big(\Sigma_{[j^*],[j^*]}\big)^{-1}_{[j^*-1],j^*}\big(\Sigma_{[j^*],[j^*]}\big)^{-1}_{j^*,[j^*-1]}}{\big(\Sigma_{[j^*],[j^*]}\big)^{-1}_{j^*,j^*}}&\geq 0,
\end{split}
\]
from which the inequalities in equation~\eqref{eqn:inductionequation}
are proven for $j=j^*$. Given that the first inequality in equation~\eqref{eqn:inductionequation}
is equivalent to the function in~\eqref{eqn:middle} being non-decreasing
in its arguments, we have proven the required result.
\end{proof}

Lemma~\ref{thm:reorderingaCIS} allows us to find a row
of the precision matrix $K$ whose off-diagonals are non-positive and
assume it is the last element of a CIS ordering. This is the basis of
our algorithm.

\begin{theorem}
Suppose $X$ is a centered multivariate Gaussian for which there exists 
a CIS ordering. Then the following procedure produces a permutation 
$\sigma$ such that $X_\sigma$ is CIS.
\begin{enumerate}
\item Initialize $O^{(1)}=[m]$ as the ``leftover'' set, $K^{(1)} = K$ as
the current precision matrix, and $C^{(1)} = \{j: K_{j,\backslash j}\leq 0\}$
as the current candidate set.
\item For $i=1,\ldots, m$, take an element of $k \in C^{(i)}$ and set $\sigma(m-i+1) = k$.
Compute
\begin{align*}
O^{(i+1)} &= O^{(i)}\backslash \{k\},\\
K^{(i+1)} &= \Big(\Sigma_{O^{(i+1)},O^{(i+1)}}\Big)^{-1},\\
C^{(i+1)} &= \{j : K^{(i+1)}_{j,\backslash j} \leq 0\}.
\end{align*}
\end{enumerate}
\end{theorem}
\begin{proof}
We must simply show that at each step $C^{(i)}$ is not empty, since at 
each step, the condition $K^{(i)}_{j,\backslash j} \leq 0$ is 
sufficient for the variable $X_j$ to be a non-decreasing function in 
the variables $X_v$ with $v \in O^{(i)}\backslash\{j\}$ by 
Lemma~\ref{lem:nonnegprecis}. This follows by existence of a CIS 
ordering along with Theorem~\ref{thm:reorderingaCIS}. Indeed if a CIS 
ordering exists, then $C^{(1)} \neq \emptyset$, in which case, an 
\textit{arbitrary} element of $C^{(1)}$ can be taken to be $\sigma(m)$. A simple induction argument shows that this is true for 
each $C^{(i)}$. 
\end{proof}

We illustrate this algorithm with an example.

\begin{example}\label{ex:4cycle0}
Consider the four dimensional Gaussian distribution with covariance and precision matrix
	$$
	\Sigma\;=\;\begin{bmatrix}
		1 & 0.75 &  0.50 &  0.14\\
		0.75 &  1 &   0.81 &  0.50\\
		0.50 &   0.81 &  1 &  0.75\\ 
  0.14 &  0.50 &  0.75 &  1
	\end{bmatrix},\qquad K=\begin{bmatrix}
		2.77 & -2.51 & 0 & 0.88\\
		-2.51 & 5.49 & -3.2 & 0\\
		0 & -3.2 & 5.49 & -2.51\\
		0.88 & 0 & -2.51 & 2.77
	\end{bmatrix}.
 	$$ 
The matrix $K$ has two rows with only non-positive off-diagonal entries. We choose $i_1=2$ and consider the marginal distribution over $\{1,3,4\}$. The matrix 
	$$(\Sigma_{134})^{-1}=\begin{bmatrix}
	1.61 &-1.47 & 0.88\\
	-1.47&3.62&-2.51\\
	0.88&-2.51&2.77
\end{bmatrix}
$$  
has one row with nonpositive off-diagonal entries; so we take $i_2=3$. This shows that both $(1,4,3,2)$ and $(4,1,3,2)$ are CIS orderings. Beginning with $i_1=3$ shows that also $(1,4,2,3)$ and $(4,1,2,3)$ are CIS orderings and there are no other CIS orderings of $X$.
\end{example}

\subsection{Noisy CIS Recovery}

In the noisy setting, we are given a matrix of observations
$\BX \in \BBR^{n\times m}$ where the rows are i.i.d 
and distributed according to $\cN_m(0,\Sigma)$,
where $\Sigma$ is such that the distribution admits
a CIS ordering. As in Section~\ref{sec:choleskestimat}, we let $\mathbf{x}_t$ refer to the $t$-th column of $\BX$. For any $i\in \{1,\ldots,m\}$ and any non-empty $A\subseteq \{1,\ldots,m\}\setminus \{i\}$, denote by $\beta^{(i,A)}$ the vector of coefficients of the linear regression  of $\mathbf{x}_i$ on $\BX_{[n],A}$. Then we have that
\[
\beta^{(i,A)}\;=\;\Sigma_{i,A}\Sigma_{A,A}^{-1}.
\]
When $\beta^{(i,A)}\geq 0$, we say $\mathbf{x}_i$ can be positively regressed on $\BX_{[n],A}$.
For $\alpha > 0$, an estimator $\hat{\beta}^{(i,A)}$  (we suppress $n$-dependence for ease) of $\beta^{(i,A)}$ is said to be
$n^\alpha$,  consistent if
\[
n^\alpha(\hat{\beta}^{(i,A)} - \beta^{(i,A)}) \to 0
\]
in probability as $n\to\infty$. 
Our noisy CIS recovery algorithm presented in the Theorem
below will mimic the method of the previous section by
inspecting the entries of $\hat{\beta}^{(i,A)}$ at each step, assuming a bound on the entries of $\beta^{(i,A)}$.

\begin{theorem}
Assume that there exists a CIS ordering of the distribution $\cN_m(0,\Sigma)$
and there exists an $\epsilon^* = \epsilon^*(\Sigma) > 0$ such that for any $i\in V$ and $A\subseteq V\setminus \{i\}$, either $\beta^{(i,A)}$ is a non-negative vector or $\min_j \beta^{(i,A)}_j<-2\epsilon^*$. For an $\alpha > 0$, let $\hat{\beta}^{(i,A)}$ be an $n^\alpha$-consistent estimators of $\beta^{(i,A)}$ and let $\epsilon_n$ be a sequence such that $\epsilon_n\to 0$ while $n^\alpha \epsilon_n \to \infty$.

We define an estimator $\hat\sigma$ through the following algorithm:
\begin{enumerate}
    \item Initialize $\mathcal{A}_1 = [m]$ to be the set of active variables and set $t=1$.
    \item If $t \leq m-2$, for each $i \in \mathcal{A}_t$, we compute $\hat{\beta}^{(i,\mathcal{A}_t\backslash\{i\})}$. At the first instance\footnote{In practice we could score different potential choices to further improve the power of the method.} of $i^*$
    such that all entries of $\hat\beta^{(i^*,\mathcal{A}_t\backslash\{i^*\})}$ are greater than $-\epsilon_n$ we define
    \[
    \hat{\sigma}(m-t+1) = i^*.
    \]
    Define $\mathcal{A}_{t+1} = \mathcal{A}_t \backslash \{i^*\}$ and increment $t$ and repeat this step until $t = m-1$.
    \item When $t=m-1$ it must be that $|\mathcal{A}_{t}|=2$, in which case, we take $\hat{\sigma}(1)$ and $\hat{\sigma}(2)$ to be arbitrary.
\end{enumerate}
As $n\to\infty$,  $\hat{\sigma}$ will be a valid CIS ordering of $\cN_m(0,\Sigma)$ with probability going to 1. 
\end{theorem}
\begin{proof}
    Depending on the sample size $n$, consider the event
\begin{equation}\label{eq:Enn}
\mathcal E^{(n)}\;=\;\bigcap_{i,A}\mathcal E_{i,A}^{(n)},\qquad \mathcal E_{i,A}^{(n)}:=\{\|\hat\beta^{(i,A)}-\beta^{(i,A)}\|_\infty<\epsilon_n\}.    
\end{equation}
By $n^\alpha$-consistency of the estimators and the fact that $n^\alpha\epsilon_n\to \infty$, $\P(\mathcal E^{(n)})\to 1$ as $n\to \infty$\footnote{Indeed if $A_n$, $B_n$ are sequences of events such that $\P(A_n)\to 1$ and $\P(B_n)\to 1$ then $\P(A_n\cap B_n)\to 1$ simply because $\P((A_n\cap B_n)^c)=\P(A_n^c\cup B_n^c)\leq \P(A_n^c)+\P( B_n^c)\to 0$.}. 
Note that, by the definition of $\epsilon^*$ and the fact that $\epsilon_n<\epsilon^*$ if $n$ is sufficiently large, conditionally on $\mathcal E^{(n)}$, this is equivalent to the fact that $\mathbf{x}_i$ can be positively regressed on $\BX_{[n],A}$. More specifically, let $R_t$ for $t=1,\ldots, m-3$, be the event that says that at the $t$-th step of the algorithm:
\begin{itemize}
    \item [(a)] $\BX_{[n],\mathcal A_t}$ admits a $\CIS$ ordering,
    \item [(b)] the algorithm correctly finds an $\mathbf{x}_i$ that can be positively regressed on $\BX_{[n],\mathcal A_t\setminus \{i\}}$.
\end{itemize}
Note that (a) is automatically satisfied if $t=1$. Similarly, for an arbitrary $t$, (a) holds automatically conditionally on $R_1\cap\ldots\cap R_{t-1}$, by Theorem~\ref{thm:reorderingaCIS}. The probability of recovering a $\CIS$ ordering is $\P(R_1\cap \ldots \cap R_{m-3})$ and we have
$$
\P(R_1\cap \ldots \cap R_{m-3})\;=\;\P(R_1)\P(R_2|R_1) \cdots \P( R_{m-3}|R_1\cap \ldots\cap R_{m-4}).
$$
Denote $\P^{(n)}(\cdot)=\P(\cdot|\mathcal E^{(n)})$. We also have
$$
\P^{(n)}(R_1\cap \ldots \cap R_{m-3})\;=\;\P^{(n)}(R_1)\P^{(n)}(R_2|R_1) \cdots \P^{(n)}( R_{m-3}|R_1\cap \ldots\cap R_{m-4}).
$$
As we said earlier, after conditioning on $\mathcal E^{(n)}$, $X_i$ can be positively regressed on $X_A$ if and only if all coefficients of $\hat \beta^{(i,A)}$ are greater than $-\epsilon$. This means that
$$
\P^{(n)}(R_1)\;=\;\P^{(n)}(R_2|R_1)\;=\;\cdots\;=\; \P^{(n)}( R_{m-3}|R_1\cap \ldots\cap R_{m-4})\;=\;1
$$
implying that $\P^{(n)}(R_1\cap \ldots \cap R_{m-3})=1$. This implies that $\P(R_1\cap \ldots \cap R_{m-3})\to 1$ as $n\to\infty$\footnote{Indeed, if $A_n, B_n$ are sequences of events such that $\P(A_n|B_n)=1$ and $\P(B_n)\to 1$ then $\P(A_n\cap B_n)=\P(B_n)\to 1$. Since $\P(A_n)\geq \P(A_n\cap B_n)$, then also $\P(A_n)\to 1$.}, which completes the proof.
\end{proof}

\begin{remark}
The event $\mathcal E^{(n)}$ in \eqref{eq:Enn} may have small 
probability for finite sample sizes. However, for the proof it is not 
necessary to define $\mathcal E^{(n)}$ as an intersection over all 
pairs $(i,A)$. For example, it is sufficient to include only the pairs 
$(i,A)$ such that $\BX_{[n],A\cup \{i\}}$ admits a $\CIS$ ordering but 
$\mathbf{x}_i$ cannot be positively regressed on $\BX_{[n],A}$ (if 
$\Sigma$ is an inverse M-matrix then there are no such pairs). 
\end{remark}

\bibliographystyle{plain}
\bibliography{biblio}

\end{document}